\documentclass[a4paper,11pt]{article}
\usepackage{amsfonts,amssymb,latexsym,amsmath}
\usepackage{multicol}
\newcounter{num}
\newcommand{\Num}{\refstepcounter{num}%
\textbf{
\arabic{num}}}

\newcommand{\Theorem}{\textbf{Theorem ~}}
\newcommand{\Prop}{\textbf{Proposition~}}
\newcommand{\Proof}{\textbf{Proof}}

\newcommand{\Lemma}{\textbf{Lemma~}}
\newcommand{\Def}{\textbf{Definition~}}
\newcommand{\Cor}{\textbf{Corollary~ }}

\newcommand{\al}{\alpha}

\newcommand{\la}{\lambda}

\newcommand{\gx}{{\mathfrak g}}
\newcommand{\ogx}{\overline{\gx}}
\newcommand{\ohx}{\overline{\hx}}
\newcommand{\zx}{{\mathfrak z}}
\newcommand{\px}{{\mathfrak p}}

\newcommand{\hx}{{\mathfrak h}}

\newcommand{\ut}{{{\mathfrak u}{\mathfrak t}}}

\newcommand{\UT}{{{\mathrm U}{\mathrm T}}}

\newcommand{\Ad}{{\mathrm{Ad}}}
\newcommand{\ad}{{\mathrm{ad}}}
\newcommand{\Ker}{{\mathrm{Ker}}}
\newcommand{\codim}{{\mathrm{codim}}}
\newcommand{\mult}{{\mathrm{mult}}}
\newcommand{\ind}{{\mathrm{ind}}}
\newcommand{\Tr}{{\mathrm{Tr}}}

\newcommand{\Fb}{{\Bbb F}}
\newcommand{\Cb}{{\Bbb C}}

\renewcommand{\leq}{\leqslant}
\renewcommand{\geq}{\geqslant}

\begin{document}
\Large

\title{The orbit method for unipotent groups over finite field}
\author{A.N.Panov
\thanks{The paper is supported by the RFBR-grants 12-01-00070,
12-01-00137, 13-01-97000-Volga region-a}}

\date{}

\maketitle

According to  A.A.Kirillov's orbit method there exists one to one
correspondence between the irreducible representations of an
arbitrary connected, simply connected Lie group and its coadjoint
orbits. This correspondence between orbits and representations makes
 possible to solve problems of the representation theory in
terms of coadjoint orbits. The orbit method initiate many  papers
beginning from 1962. It turns out that the ideas of the orbit method
are useful for  the large classes of Lie groups  (see \cite{V1,V2}),
and also for some matrix groups, defined over finite field.

In our paper we obtain formula for multiplicities of certain
representations of unipotent groups over finite field in terms of
coadjoint orbits (see theorem  \ref{tt} and corollaries). For
reader's convenience we formulate and prove the maim statements of
the orbit method over finite field (see \cite{Kzh}).

 Let  $K =\Fb_q$ be a finite field of characteristic $p$ having $q=p^m$ elements.
Let $\gx$ be a  subalgebra of the Lie algebra $\ut(N,K)$, consisting
of all upper triangular matrices  with zeros on the diagonal.
Suppose that $p$ is large enough to determine the exponential
$\exp(x)$ map on $\gx$. For instance, let $p\geq N$. Then the
exponential map is a bijection of the Lie algebra  $\gx$ onto the
subgroup $G=\exp(\gx)$ of the unitriangular group $\UT(N,K)$. One
can define the adjoint representation  of the group $G$ on
 $\gx$ by the formula  $\Ad_g(x) = gxg^{-1}$.

 Denote by  $\gx^*$ the conjugate space of $\gx$. One can define the
coadjoint representation  of the group $G$ in $\gx^*$ by the formula
$\Ad^*_g\la(x) = \la(\Ad_g^{-1} x)$.

 Note that if  $\gx^\la$ is a stabilazer of  $\la\in\gx^*$,
 then  the subgroup  $G^\la=\exp(\gx^\la)$ is a stabilizer of $\la$
 in  $G$. One can calculate the number of elements $|\Omega|$ of the orbit  $\Omega=\Ad_G^*(\la)$ by the formula
 \begin{equation}
|\Omega| = \frac{|G|}{|G^\la|} = \frac{|\gx|}{|\gx^\la|} = q^{\dim
\gx - \dim \gx^\la} = q^{\frac{1}{2}\dim\Omega}.
 \end{equation}
\Def\Num. A subalgebra  $\px$ of $\gx$ is a polarization of $\la\in
\gx^*$, if $\px$ is a maximal isotropic subspace for the skew
symmetric bilinear form  $B_\la(x,y) = \la([x,y])$ on $\gx$. Recall
that the subspace $\px$ is isotropic, if  $B_\la(x,y) =0$ for any
$x,y\in\px$.

Note that any polarization contains the stabilizer  $\gx^\la$,
because  $\px +\gx^\la$ is an isotropic subspace.\\
\Prop\Num. Any linear form $\la$ on a nilpotent Lie algebra  $\gx$ has a polarization.\\
 \Proof. We shall prove using induction method for the dimension of the Lie algebra  $\gx$. The statement if obvious for
 one dimensional Lie algebras, since in this case  Lie algebra is a polarization.
 Assume that the statement is proved for  all Lie algebras of dimension  $<\dim(\gx)$. We are going to prove  the statement
 for  $\dim(\gx)$.

 If the dimension of a center  $\zx$ the Lie algebra  $\gx$ greater than one, then
 one can prove existence of polarization applying induction assumption for the factor algebra
 of   $\gx$ with respect to the ideal  $\Ker(\la|_\zx)$.
 Similarly, for the case  $\dim(\zx)=1$,~~ $\la|_\zx
 =  0$.

 Let  $\zx =Kz$ and $\la(z)\ne 0$. Consider the two dimensional ideal  $Ky+Kz$,
 containing  $\zx$. There exists a character  $\al$ of the Lie algebra $\gx$ such that
  $\ad_u(y) = \al(u)z$ for any  $u\in\gx$. The kernel $\gx_0$ of the character
  $\al$ is an ideal of codimension one in  $\gx$. There exists an element
  $x\in\gx$ such that $\gx= Kx+\gx_0$ and $[x,y]=z$.

  Denote by   $\la_0$ the restriction of  $\la$ on $\gx_0$.
  According the induction assumption  $\la_0$ has a polarization $\px_0$ in $\gx_0$.
  Let us prove that  $\px_0$ is also a polarization for  $\la$ in $\gx$. Really,
  $\px_0$ is a subalgebra and a maximal isotropic subspace in   $\gx_0$; we will show that $\px_0$ is a maximal isotropic subspace
  in  $\gx$.
  Suppose that one can extent  $\px_0$ to  an isotropic subspace adding the element  $x+u_0$, where  $u_0\in\gx_0$. Note that  $z,y$
   belong to the stabilizer  $\gx_0^{\la_0}\subset \px_0$.
   Then  $0 =\la([x+u_0,y]) = \la([x,y]) = \la(z) \ne 0$. A contradiction.
   $\Box$\\
   \Prop\Num\label{pp}. Let  $\px$ be a polarization of $\la\in\gx^*$,
   ~ $P=\exp(\px)$,~~
    $\Omega(\la)$  be the coadjoint orbit of  $\la$, ~$\pi$ be the natural projection of
   $\gx^*$ onto $\px^*$, ~ $L^\la = \pi^{-1}\pi(\la)$. Then \\
   1)~ $\dim \px = \frac{1}{2}\left( \dim \gx + \dim\gx^\la\right)$;   \\
   2) ~$|L^\la| = \sqrt{|\Omega(\la)|}$;\\
   2)~ $L^\la = \Ad^*_P \la$, in particular  $ L^\la\subset  \Omega(\la)$.\\
   \Proof. The statement 1) follows from the formula of dimension of a
   maximal isotropic subspace for the skew symmetric bilinear form
    $B_\la(x,y)$. From  1) we obtain
   $$\codim\,\px =  \frac{1}{2}\left( \dim \gx
   - \dim \gx^\la\right) = \frac{1}{2}\dim \Omega(\la).$$
   This implies the statement  2):
   $$|L^\la| = q^{\codim\,\px} = q^{\frac{1}{2}\dim
   \Omega(\la)} = \sqrt{|\Omega(\la)|}. $$

   Since  $\la([x,y])=0$ for any  $x,y\in\px$, we have  $\ad_\px^* \la (y)
   = 0$ for any  $y\in\px$. Then  $\Ad^*_P\la(y) = \la(y)$  for any  $y\in \px$.
   This is equivalent to $$ \Ad^*_P\la \subset L^\la. $$

   The equality   $ \Ad^*_P\la = L^\la $ is true, since this subsets have equal number of elements:
   $$
   |\Ad^*_P\la| = \frac{|P|}{|G^\la|} = q^{\dim \px - \dim \gx^\la} =  q^{\frac{1}{2}\dim
   \Omega(\la)} = |L^\la|.\,
   \Box $$

   Fix a non trivial  character
   $e^x: K\to\Cb^*$. We have    \begin{equation}\label{exp}
   \sum_{t\in\Fb_q} e^{\al t} = \left\{\begin{array}{l}
   q,~\mbox{если}~~ \al=0,\\
0,~\mbox{если}~~ \al\ne 0.\end{array}\right.\end{equation} The
equality  (\ref{exp}) is easy to prove: the image of homomorphism
$e^x$ is a subgroup  of $\Cb^*$; if  $\al\ne 0$, then this subgroup
is nontrivial and coincides with  the subgroup of all roots of some
order $m\ne 1$ of unity; the sum of all roots of order  $m\ne 1$ of
unity equals to zero.

   Restriction of  $\la$ on its polarization  $\px$
   defines a character (one dimensional representation)  $\xi$ of the group
   $P=\exp(\px)$ by the formula
   $$\xi_\la(\exp(x)) = e^{\la(x)}.$$
Consider the induced representation
\begin{equation}
T^\la = \ind (\xi_\la, P, G).
\end{equation}
Denote by  $\chi_\la(g) = \Tr~ T^\la(g)$ the character of
representation
$T^\la$.\\
\Theorem\Num\label{chiq}. $$\chi_\la(g) = \frac{1}{\sqrt{|\Omega|}}
\sum_{\mu\in\Omega(\la)} e^{\mu(\ln(g))}$$ \Proof. Extent $\xi_\la$
from $P$ to $G$ by the formula
$$\tilde{\xi}_\la(u) = \left\{ \begin{array}{cl} \xi_\la(u)&, ~~ \mbox{if}
~~ u\in P,\\
0&, ~ ~\mbox{if} ~~ u\notin P.\end{array}\right.$$ Formula
(\ref{exp}) implies
$$ \tilde{\xi}_\la(u) = \frac{1}{|L^\la|}\sum_{\mu\in L^\la}
\xi_\mu(u) = \frac{1}{|L^\la|}\sum_{p\in P} \xi_{\Ad^*_p(\la)}(u).$$

Choose the system of representatives $\{g_i:~ i=\overline{1, k}\}$
of the classes  $G/P$. Using the well known formula for induced
characters (see \cite[chapter 6]{KR}), we obtain
$$\chi_\la(u) = \sum_{g^{-1}_iug_i\in P}\tilde{\xi}_\la(g^{-1}_iug_i) =
\frac{1}{|L^\la|} \sum_{i=1, p\in P}^k \xi_{\Ad_p^*\la}
(g^{-1}_iug_i) = \frac{1}{|L^\la|} \sum_{i=1, p\in P}^k
\xi_{\Ad_{g_ip}^*\la}(u),
$$
Finally,
$$ \chi_\la(u) = \frac{1}{|L^\la|} \sum_{g\in G} \xi_{\Ad_g^*\la}(u)
= \frac{1}{\sqrt{|\Omega|}} \sum_{\mu\in\Omega(\la)}
e^{\mu(\ln(u))}. \Box
$$
\Theorem\Num\label{t}. \\
1) $ \dim T^\la = q^{\frac{1}{2}\dim \Omega(\la)} = \sqrt{|\Omega|}$.\\
2) The representation  $T^\la$ does not depend on the choice of polarization.\\
3) The representation $T^\la$ is irreducible.\\
4) Representations  $T^\la$ and $T^{\la'}$ are equivalent if and
only if  $\la$ and $\la'$ belong to the same  coadjoint orbit.\\
5) For any irreducible representation $T$ of the group  $G$ there
exists  $\la\in\gx^*$ such that the representation  $T$ is equivalent to
 $T^\la$.\\
\Proof. The statement  1) follows from  $$ \dim T^\la = \dim
(\ind(\xi_\la, P,G)) = q^{\codim\,\px}.$$ The statement  2) is a
corollary of  the theorem \ref{chiq}.

Let us show that the system of characters  $\{\chi_\la\}$, where
$\la$ is running through some  system of representatives of the
coadjoint orbits, is orthonormal. Let  $\Omega$, ~$\Omega'$ be two
coadjoint orbits and
 $\la$, ~$\la'$ be representatives of this orbits. Then
$$
(\chi_\la, \chi_{\la'}) = \frac{1}{|G|}\sum_{u\in
G}\chi_\la(u)\overline{\chi_{\la'}(u)}= \frac{1}{|G|}\cdot
\frac{1}{\sqrt{|\Omega|\cdot|\Omega'|}}\cdot \sum_{\mu\in\Omega,
\mu'\in\Omega', u\in G}\xi_\mu(u)\overline{\xi_{\mu'}(u)}=
$$
$$
\frac{1}{|G|}\cdot \frac{1}{\sqrt{|\Omega|\cdot|\Omega'|}}\cdot
\sum_{\mu\in\Omega, \mu'\in\Omega', u\in G}e^{(\mu-\mu')\ln(u)}.
$$
Applying  $$ \sum_{x\in\gx}e^{\eta(x)}=\left\{\begin{array}{cl}
|G|&, ~~\mbox{if}~~ \eta=0,\\
0&, ~~\mbox{if}~~ \eta\ne 0 \end{array}\right.,$$ we obtain that, if
$\Omega \ne \Omega'$, then $(\chi_\la, \chi_{\la'})=0$.

In the case  $\Omega=\Omega'$, we have got
$$ (\chi_\la, \chi_\la) = \frac{1}{|G|\cdot |\Omega|}
\sum_{\mu, \mu'\in\Omega}\sum_{x\in\gx}e^{(\mu-\mu')x} =
\frac{1}{|G|\cdot |\Omega|}\cdot  |\Omega| \cdot |G| = 1.
$$
This proves  3) and 4).

We shall use notation  $T^\Omega$ for the class of equivalent
representations  $T^\la$, where $\la\in\Omega$. To prove  statement
5) we verify that the  sum of squares of dimensions of irreducible
representations  $\{ T^\Omega: ~ \Omega\in\gx^*/G\}$ equals to the
number of elements of the group:\\
$$\sum_{\Omega\in\gx^*/G} \left(\dim T^\Omega\right)^2 =
\sum_{\Omega\in\gx^*/G} \left( \sqrt{|\Omega|} \right)^2 =
\sum_{\Omega\in\gx^*/G}|\Omega| = |\gx^*| = |G|.
$$
This proves 5).$\Box$

\Lemma\Num\label{lsa}. Let  $\gx$ be a nilpotent Lie algebra,
$\gx_0$  be a subalgebra of  $\gx$ of codimension one. Then $\gx_0$
is an ideal of $\gx$.
\\
\Proof. Suppose the contrary. Then  $[\gx,\gx_0]\ne\gx_0$; the exist
the elements  $y\in\gx_0$, ~$x\notin\gx_0$ such that  $[x,y] = \al x
\bmod\gx_0$,~ $\al \ne 0$. The subalgebra $\gx_0$ is invariant with
respect to  $\ad_{y_0}$. Since $\gx = kx\oplus\gx_0$, the operator
$\ad_{y_0}$ is not nilpotent in $\gx/\gx_0$; this contradicts to
assumption that the Lie algebra $\gx$ is nilpotent. $\Box$

\Lemma\Num\label{lto}. Let  $\gx$,~$\gx_0$ be as in lemma
\ref{lsa},~ $\pi$ is a projection  $\gx^*\to\gx_0^*$; ~
$\la_0\in\gx_0^*$, ~ $\omega = \Ad_{G_0}^*(\la_0)$, ~ $\gx^{\la_0}
=\{ x\in\gx:~ \la_0([x,\gx_0)] = 0\}$.\\
 1) Let the subalgebra  $ \gx^{\la_0}$ belong to  $\gx_0$. Then\\
 1a)~ $\pi^{-1}(\la_0)$ lie in the same содержится  $\Ad^*_G$-orbit
 $\Omega$;\\
 1b) ~  $\dim \Omega = \dim \omega +2$ (i.e. $|\Omega| =
 q^2|\omega|$) .\\
 2) Let the subalgebra  $ \gx^{\la_0}$  do not lie in $\gx_0$. Then for any
  $\Ad_G^*$-orbit $\Omega$, which has nonempty intersection with
   $\pi^{-1}(\la_0)$, the projection  $\pi$ establishes one to one correspondence
   between  $\Omega$  and $\omega$; in particular,  $\dim \Omega = \dim \omega$.\\
 \Proof.\\
 1) Suppose that the subalgebra  $ \gx^{\la_0}$ belongs to $\gx_0$.
 Since  $[\gx,\gx]\subset\gx_0$, the formula $B_0(x,y) = \la_0([x,y])$ defines
 a skew symmetric bilinear form on the Lie algebra  $\gx$.
The kernel  $V$ of the bilinear form  $B_0$ coincides with
$\gx^{\la_0}$ and belongs to  $\gx_0$. The kernel  $V_0$ of the
restriction  $B_0$ on $\gx_0$ coincides with $\gx_0^{\la_0}$.

 Let us prove that  there exists a pair of elements
    $u\in \gx\setminus\gx_0$  and $v\in V_0$,
 such that   $B_0(u,v)=1$. Really, decompose  $\gx_0 =
 L_0\oplus V_0$, where  $L_0$ is a subspace, with the property that the restriction
 of bilinear form  $B_0$ on $L_0$ is nondegenerate. Choose an arbitrary element
 $u'\in\gx\setminus \gx_0$ and consider the linear form
$B_0(u',\cdot)$ on  $L_0$. There exists  $x_0\in L_0$ such that $
B_0(u',\cdot) = B_0(x_0,\cdot)$ on $L_0$. The element  $u=u'-x_0$
satisfies $B_0(u,L_0)=0$. Since $u\notin V$, there exists  $v\in
V_0$ such that $B_0(u,v)=1$.

By direct calculations, we verify that for  any
$\la\in\pi^{-1}(\la_0)$ the following equalities are valid
\begin{equation}\label{first}
\left\{\begin{array}{l}\Ad^*_{\exp(tv)}\la(u) = \la(u)+t,\\
 \Ad^*_{\exp(tv)}\la(y) = \la(y)
~~\mbox{for~~ any}~~~ y\in\gx_0.
\end{array}\right.\end{equation}
This implies the statement 1a).

The orbit  $\Omega$ is a union
\begin{equation}\label{omegaf}\Omega = \bigcup_{t\in K}\pi^{-1}(\omega_t),
\end{equation}
where $\omega_t = \Ad^*_{\exp(tu)}\omega$  is a coadjoint orbit in
$\gx_0^*$. Let us show that the orbits $\omega_t$ are pairwise
different. Really, if not, there exists  $t'\ne t''\in K$ and
$g'_0,g_0''\in G_0$ such that
$$\Ad^*_{\exp(t'u)}\Ad^*_{g_0'}\la_0 =
\Ad^*_{\exp(t''u)}\Ad^*_{g''_0}\la_0.$$ Then
$$ \Ad^*_{\exp(tu)}\Ad^*_{g_0}\la_0 = \la_0,$$
where $t=t'-t''\in K^*$ and  $g_0$ is an  element of $G_0$. Then the
stabilizer  $\exp(\gx^{\la_0})$ does not belong to $G_0$; this
contradicts to  assumption of the item  1). Using  (\ref{omegaf}),
we obtain $|\Omega| = q^2|\omega|$. Therefore $\dim \Omega = \dim
\omega + 2$. This proves   1b).

Turn to proof of the item  2). Suppose that the subalgebra $
\gx^{\la_0}$ does not belong to   $\gx_0$.  For an arbitrary nonzero
element  $x$ of $ \gx^{\la_0}\setminus \gx_0$ we have decomposition
$\gx = Kx\oplus\gx_0$. The group  $G$ is a semidirect product
$G=G_0X$, where $X=\{ \exp(tx):~ t\in K\}$.

 Let $\la\in\pi^{-1}(\la_0)$. The equality
$$
\la([x,\gx]) = \la([x,Kx\oplus\gx_0]) = \la_0([x,\gx_0]) =0$$
implies that  $x$ belongs to  $\gx^\la$. The subgroup  $X$ lies in
the stabilizer of  $\la$. The orbit  $\Omega(\la)$ coincides with
$\Ad^*_{G_0}(\la)$. Since the projection $\pi:\gx^*\to\gx_0^*$ is
invariant with respect to  $\Ad^*_{G_0}$, the map  $\pi$ project
$\Omega(\la)$ onto $\omega$.

It remains to show that
\begin{equation}\label{ee} \Omega(\la) \cap\pi^{-1}(\la_0) =
\{\la\}\end{equation} for any  $\la\in\pi^{1}(\la_0)$. Suppose that
 $\la'=\Ad^*_g\la$ and $\la,\la'\in\pi^{-1}(\la_0)$. Since  $G
= G_0X$ and $X\in G^\la$, we verify  that $\la' = \Ad^*_{g_0}\la$
for some  $g_0\in G_0$. As $\la,\la'\in\pi^{-1}(\la_0)$, the element
$g_0$ lies in stabilizer  $G_0^{\la_0}$. Then $g_0 = \exp(y_0)$ for
some  $y_0\in \gx_0^{\la_0}$. We obtain
\begin{equation} \label{ll} \la'(x) = \la(\Ad^*_{\exp(-y_0)}x) = \la(x)
- \la(\ad_{y_0}x) + \sum_{k\geq
2}\frac{(-1)^k}{k!}\la_0(\ad_{y_0}^kx).
\end{equation}
Since  $x\in \gx^\la$,  we have $ \la(\ad_{y_0}x) =0$. As $y_0\in
\gx_0^{\la_0}$, we have $\la_0(\ad_{y_0}^kx)=0$ for any $k \geq 2$.
Substituting into  (\ref{ll}), we obtain  $\la'(x) = \la(x)$. Using
 $\la,
\la'\in\pi^{-1}(\la_0)$, we conclude  $\la = \la'$. $\Box$\\
\Lemma\Num\label{ss}. For any subalgebra  $\hx$ of a nilpotent Lie
algebra   $\gx$ there exists a chain of subalgebras  $\gx=\gx_0
\supset\gx_1\supset\ldots\supset\gx_k=\hx$ such that  $\gx_{i+1}$ is
an ideal of codimension one in  $\gx_{i}$ for any  $1\leq i\leq
k-1$.
\\
\Proof. We use the induction method for  $\dim \gx$. For  $\dim \gx
=1$ the statement is obvious. Assume  that the statement is true for
$\dim \gx = n-1$; let us prove it  for  $n$. The nilpotent Lie
algebra $\gx$ has a nonzero central element  $z$. Consider
 projection  $\phi:\gx\to\ogx = \gx/Kz$. The image  $\ohx$ is a subalgebra in $\ogx$.
 As $\dim\ogx <n$, according to induction assumption, there exists a chain of subalgebras
  $\ogx=\ogx_0
\supset\ogx_1\supset\ldots\supset\ogx_k=\ohx$, where $\ogx_i$ is an
ideal of  codimension one in $\ogx_{i+1}$ for any  $1\leq i\leq
k-1$. Denote  $\gx_i = \phi^{-1}(\ogx_i)$. If $z\in\hx$, then $\gx_k
= \hx$; this completes  construction of the chain of subalgebras.
 If $z\notin\hx$, then  $\gx_k = \hx+Kz$. It remains to put
$\gx_{k+1}$ equal to  $\hx$. $\Box$
 \\
  \Theorem\Num\label{tt}. Let  $G=\exp(\gx)$ be an unipotent group over the finite field
  $K$, ~  $\hx$ be a subalgebra of $\gx$, ~  $H=exp(\hx)$.
Let  $\Omega$ (resp. $\omega$) be a coadjoint orbit in $\gx^*$
(resp. $\hx^*$), ~$T^\Omega$ and $t^\omega$ be  corresponding
irreducible representations  of $G$ and $H$, $\pi$ be the natural
projection $\gx^*$ onto  $\hx^*$. Denote  $m(\omega, \Omega) =
\mathrm{mult} (T^\Omega,
 \mathrm{ind}(t^\omega,G)) = \mathrm{mult} (t^\omega,
 \mathrm{res}(T^\Omega,H))$. Then
   $$m(\omega, \Omega)
 =\frac{|\pi^{-1}(\omega)\cap\Omega|}{\sqrt{|\omega|\cdot |\Omega|}}.$$
  \Proof. Introduce notations  $$ P =
|\pi^{-1}(\omega)\cap \Omega|,\quad\quad  Q = \sqrt{|\omega|\cdot
|\Omega|},\quad\quad M = \mult\left(T^\Omega, \ind(t^\omega,
G)\right).$$

We shall prove  that $M=P/Q$ using the induction  method with
respect to $\codim(\hx,\gx)$. If  $\codim(\hx,\gx)=0$, then  $\Omega
= \omega$ and hence $P = |\Omega|$, ~$Q = |\Omega|$, ~$M=1$; this
proves the equality  $M=P/Q$.

Assume that the equality is proved for $\codim(\hx,\gx)< k$; let us
prove for  $\codim(\hx,\gx)=k$. The lemma  \ref{ss} implies that
there exists a  subalgebra   $\gx_1$ obeying the conditions
$\gx\supset\gx_1\supset\hx$,~ $\codim(\hx,\gx_1)=1$. Choose
$\la_0\in\omega$. The natural projections  $\pi:\gx^*\to\hx^*$,~
$\pi_1:\gx^*\to\gx_1^*$, ~$\pi_0:\gx_1^*\to\hx^*$ satisfy $\pi =
\pi_0\pi_1$. For $\gx_1^{\la_0} =\{ x\in\gx_1:~ \la_0[x,\hx] = 0 \}$
only two cases are  possible: ~ $\gx_1^{\la_0}\subset \hx$, or
$\gx_1^{\la_0}\not\subset \hx$.\\
1) Case $\gx_1^{\la_0}\subset \hx$. Following  lemma  \ref{lto},
$\pi_0^{-1}(\omega)$ belongs to the same coadjoint orbit
 $\Omega_1\subset \gx_1^*$ of the group  $G_1 =\exp(\gx_1)$.
 Therefore,
$\dim \Omega_1 = \dim \omega  +2$ and
\begin{equation}\label{dd}
|\Omega_1| = q^2 |\omega| .\end{equation}

 A polarization
$\px_0$ for  $\la_0$ in $\hx$ is also a polarization for any
$\la_1\in\pi_0^{-1}$ in $\gx_1$. Really, ~$\px_0$ is an isotropic
subspace in $\gx_1$ and
$$\codim(\px_0,\gx_1) = \codim(\px_0,\hx) +1 =\frac{1}{2}(\dim \omega
+2) =\frac{1}{2} \dim \Omega_1. $$ The induced representation
$\ind(t^\omega, G_1)$ is irreducible and coincides with
$T^{\Omega_1}$,
\begin{equation}\label{ii} \ind(T^{\Omega_1},G) =
\ind(t^{\omega},G).\end{equation} According to the induction
assumption,
\begin{equation}\label{mff}
\mult\left(T^\Omega, \ind(T^{\Omega_1}, G)\right) = \frac{
|\pi^{-1}(\Omega_1)\cap \Omega|}{\sqrt{|\Omega_1|\cdot |\Omega|}}.
\end{equation}
Applying the formula (\ref{ii}), we obtain

\begin{equation} M = \frac{
|\pi^{-1}(\Omega_1)\cap \Omega|}{\sqrt{|\Omega_1|\cdot |\Omega|}}.
\end{equation}
Using  (\ref{dd}), we conclude

\begin{equation} M = \frac{ q P}{\sqrt{q^2|\omega|\cdot |\Omega|}} =\frac{P}{Q}.
\end{equation}
 2) Case  $\gx_1^{\la_0}\not\subset \hx$. By the formula  (\ref{ee})
we obtain
$$
\pi^{-1}(\omega) = \bigcup_{\la_1\in\pi^{-1}_0(\la_0)}
\pi_1^{-1}(\Omega_1(\la_1)),$$  where   $\Omega_1(\la_1)$ is an
orbit of  $\la_1\in\gx^*$ with respect to  $\Ad^*_{G_1}$. Appling
the induction assumption, we obtain
$$
P = |\pi^{-1}(\omega)\cap \Omega| =
\sum_{\la_1\in\pi^{-1}_0(\la_0)}|\pi^{-1}(\Omega_1(\la_1))\cap
\Omega| = $$ $$ \sum_{\la_1\in\pi^{-1}_0(\la_0)}
\sqrt{|\Omega_1(\la_1)|\cdot|\Omega|} \quad \mult(T^\Omega,
\ind(T^{\Omega_1(\la_1)},G)).$$ Since $|\Omega_1(\la_1)| =
|\omega|$, we have
\begin{equation}\label{PPP}
P = Q \sum_{\la_1\in\pi^{-1}_0(\la_0)}\mult(T^\Omega,
\ind(T^{\Omega_1(\la_1)},G)).\end{equation} From the other hand,
 for any polarization  $\px_0$ of  $\la_0$ the subalgebra
$$\px =\gx_1^{\la_0} +\px_0$$ is a polarization for any
$\la_1\in\pi^{-1}_0(\la_0)$.  Representation
$\ind(\xi_{\la_0},P_0,P)$ is a direct sum of one dimensional
representations  $\xi_{\la_1}$, where $\la_1\in\pi^{-1}_0(\la_0)$.
Therefore,  $ \ind(t^\omega,G)$ is a direct sum of representations
$\ind(T^{\Omega_1(\la_1)},G)$ where  $\la_1\in\pi^{-1}_0(\la_0)$. We
obtain
\begin{equation}\label{mmm} M =  \mult\left(T^\Omega, \ind(t^\omega,
G)\right) = \sum_{\la_1\in\pi^{-1}_0(\la_0)}\mult(T^\Omega,
\ind(T^{\Omega_1(\la_1)},G)).\end{equation} Substituting (\ref{mmm})
in (\ref{PPP}), we verify  $P=QM$. $\Box$\\
 \Cor\Num. The irreducible representation $T^\Omega$ occurs in decomposition  of
 $\ind(t^\omega, G)$
if and only if the orbit $\Omega$ has an nonempty intersection with
  $\pi^{-1}(\omega)$. \\
\Cor\Num. The irreducible representation  $t^\omega$ occurs in
decomposition of the restriction of representation  $T^\Omega$ on
the subgroup  $H$ if and only  if  the orbit
$\omega$ lies in   $\pi(\Omega)$. \\
\Cor\Num. Let  $\Omega$, ~$\Omega_1$,~ $\Omega_2$ be  coadjoint
orbits in $\gx^*$. Denote by  $|M|$ the number of elements in the
subset
$$ M = \{(\la_1,\la_2):~~ \la_1\in\Omega_1,~ \la_2\in\Omega_2,~
\la_1+\la_2\in\Omega \}.$$ Then
$$\mult(T^\Omega, T^{\Omega_1}\otimes T^{\Omega_2})  = \frac{
|M|}{\sqrt{|\Omega|\cdot |\Omega_1|\cdot |\Omega_2|}}.$$

{\bf Proof}. We apply theorem   \ref{tt} to the group  $G\times G$,
its coadjoint orbit  $\Omega_1\times \Omega_2$, subgroup
$H=\{(g,g):~~ g\in G\}$ and its orbit  $\omega =\{(\la,\la):~~
\la\in\Omega\}$. $\Box$

\end{document}